\documentclass{svjour3}                     % onecolumn (standard format)
\smartqed  % flush right qed marks, e.g. at end of proof
\usepackage[CJKbookmarks=true,
            bookmarksnumbered=true,
            bookmarksopen=true,
            colorlinks=true,
            citecolor=blue,
            linkcolor=blue,
            anchorcolor=green,
            urlcolor=blue
            ]{hyperref}
\usepackage{hyperref}
\usepackage{color}

\usepackage{threeparttable}
\usepackage{fmtcount}  % for footnote in table
\usepackage{amssymb}
\usepackage{amsmath}
\usepackage{booktabs}
\usepackage{graphicx}
\usepackage[latin1]{inputenc}

\usepackage[numbers]{natbib}

\newcommand{\ud}{\mathrm{d}}

\journalname{J Math Chem}
\begin{document}

\title{Numerical meshless solution of high-dimensional sine-Gordon equations via Fourier HDMR-HC approximation
}
%\subtitle{}

\titlerunning{Meshless solution of HD-SGEs via Fourier HDMR-HC}        % if too long for running head

\author{Xin Xu         \and
        Xiaopeng Luo    \and
        Herschel Rabitz %etc.
}


\institute{X. Xu \at
              Department of Chemistry, Princeton University, Princeton, NJ 08544, USA \\
              School of Managment and Engineering, Nanjing University, Nanjing, 210008, China \\
              \email{xuxin103@163.com,xx2@princeton.edu}           %  \\
           \and
           X. Luo \at
              Department of Chemistry, Princeton University, Princeton, NJ 08544, USA \\
              School of Managment and Engineering, Nanjing University, Nanjing, 210008, China\\
              \email{luo\_works@163.com,xiaopeng@princeton.edu}
           \and
           H. Rabitz \at
              Department of Chemistry, Princeton University, Princeton, NJ 08544, USA \\
              Program in Applied and Computational Mathematics, Princeton University, Princeton, NJ 08544, USA\\
              \email{hrabitz@princeton.edu}
}

\date{}

\maketitle

\begin{abstract}
  In this paper, an implicit time stepping meshless scheme is proposed to find the numerical solution of high-dimensional sine-Gordon equations (SGEs) by combining the high dimensional model representation (HDMR) and the Fourier hyperbolic cross (HC) approximation. To ensure the sparseness of the relevant coefficient matrices of the implicit time stepping scheme, the whole domain is first divided into a set of subdomains, and the relevant derivatives in high-dimension can be separately approximated by the Fourier HDMR-HC approximation in each subdomain. The proposed method allows for stable large time-steps and a relatively small number of nodes with satisfactory accuracy. The numerical examples show that the proposed method is very attractive for simulating the high-dimensional SGEs.

\end{abstract}

\section{Introduction}
\label{sec:hdmrhc.intr}
The sine-Gordon equation (SGE) is a nonlinear hyperbolic partial differential equation (PDE) involving the d'Alembert operator and the sine of the unknown function, and the SGE plays an important role in many mathematical physics applications. It was originally introduced by \citet{BourE_1862_SGE} and rediscovered by \citet{FrenkelJ&KontorovaT_1939_OnTheTheoryOfPlasticDeformationAndTwinning}. Further details about the background and applications of the SGE can be found in \cite{AeroEL&BulyginAN&PavlovYV_2009_SolutionsOf3DSGE, DrazinPG&JohnsonRS_1989_SolitonsAnIntroduction, ScottAC&ChuFYF&McLaughlinDW_1973_TheSolitonANewConceptInAppliedScience, TaleeiA&DehghanM_2014_APseudoSpectralMethodThatUsesAnOverlappingMultidomainTechniqueForNumericalSolutionOfSGEIn1and2D}. In this paper, we focus on developing an effective means for the numerical solution of the SGE in an arbitrary number of dimensions. A $(n+1)$-dimensional SGE generally takes the form:
\begin{equation}\label{eq:hdmrhc.sge}
 \begin{aligned}
  & u_{tt}(\boldsymbol{x},t)+\beta u_{t}(\boldsymbol{x},t) = \Delta u(\boldsymbol{x},t) - \psi(\boldsymbol{x})\sin(u(\boldsymbol{x},t)), \\
  &~~\boldsymbol{x}=(x_1,x_2,\cdots,x_n) \in \Omega,~t>0,
  \end{aligned}
\end{equation}
where $n$ is a positive integer, $\Delta$ is the Laplacian operator in $n$ spatial dimensions, $\Omega = [a_1,b_1]\times[a_2,b_2]\times\cdots\times[a_n,b_n] \subseteq \mathbb{R}^n$. The initial conditions associated with Eq.\eqref{eq:hdmrhc.sge} are given by
{\setlength\arraycolsep{2pt}
\begin{eqnarray}
  & u(\boldsymbol{x},0) = v_1(\boldsymbol{x}),&\boldsymbol{x}\in \Omega \label{eq:hdmrhc.sgeic1} \\
  & u_t(\boldsymbol{x},0) = v_2(\boldsymbol{x}),&\boldsymbol{x}\in \Omega, \label{eq:hdmrhc.sgeic2}
\end{eqnarray}}
and the Neumann boundary conditions are
\begin{equation}\label{eq:hdmrhc.sgebc}
 \frac{\partial u}{\partial \boldsymbol{l}}(\boldsymbol{x},t) = w(\boldsymbol{x},t),~\boldsymbol{x} \in \Gamma,~t>0
\end{equation}
where $\boldsymbol{l}$ denotes the (typically exterior) normal to the boundary of the domain, and $\Gamma$ is the boundary of $\Omega$, i.e $\Gamma = \partial \Omega$. The real parameter $\beta \geqslant 0$ weights the dissipative term. When $\beta = 0$, Eq.\eqref{eq:hdmrhc.sge} reduces to an undamped SGE in $n$ spatial variables, while when $\beta > 0$, the damped SGE is obtained. The function $\psi(\boldsymbol{x})$ can be interpreted as a Josephson current density, while $v_1$ and $v_2$ in Eqs.\eqref{eq:hdmrhc.sgeic1} and \eqref{eq:hdmrhc.sgeic2} represent wave modes or the kink and velocity, respectively.

The $(1+1)$-dimensional SGE first appeared in a strictly mathematical context in differential geometry regarding the theory of surfaces of constant curvature \cite{LambGL_1971_AnalyticalDescriptionsOfUltrashortOpticalPulsePropagationInAResonantMedium}. Moreover, it is well known that the $(1+1)$-dimensional SGE arises in many important systems, such as the Thirring model, the Coulomb gas system and the ferromagnetic $XY$ model, etc. \cite{MinnhagenP_1985_NonuniversalJumpsAndKosterlitzThoulessTransition, MinnhagenP_1987_The2DCoulombGasVortexUnbindingAndSuperfluidSuperconductingFilms, NiGJ&LouSY&ChenSQetc_1990_2DCoulombGasStudiedInSGE}. Because of its wide applications, the $(1+1)$-dimensional SGE has been studied with a variety of numerical methods, including finite difference methods (FDM) and finite element methods (FEM), etc. \cite{ArgyrisJ&HaaseM_1987_AnEngineersGuideToSolitonPhenomena}. Recently, additional solution methods have been proposed including collocation \cite{DehghanM&ShokriA_2008_ANumericalMethodFor1DNonlinearSGECollocationAndRBF, LakestaniM&DehghanM_2010_CollocationAndFiniteDifferenceCollocationMethodsForTheSolutionOfNonlinearKleinGordonEquation}, the boundary integral approach \cite{DehghanM&MirzaeiD_2008_TheBoundaryIntegralEquationApproachfForNumericalSolutionOf1DSGE, DehghanM&GhesmatiA_2010_ApplicationOfTheDualReciprocityBoundaryIntegralEquationTechniqueToSolveNonlinearKleinGordonEquation}, and a combination of the finite difference with the the diagonally implicit Runge-Kutta-Nystr\"om (DIRKN) method \cite{MohebbiA&DehghanM_2010_HighOrderSolutionOf1DSGEUsingCompactFiniteDifferenceAndDIRKNMethods}, etc.

There is recent interest in the SGE in higher dimensions. As \citet{BaroneA&EspositoF&MageeCJetc_1971_TheoryAndApplicationsOfSGE} pointed out, Eq.\eqref{eq:hdmrhc.sge} also has been applied in many branches of physics for the $n=2$ and $3$ cases. The exact solutions for the  undamped SGE in higher dimensions have been obtained by Hirota's method \cite{HirotaR_1973_ExactThreeSolitonSolutionOf2DSGE}, Lamb's method \cite{ZagrodzinskyJ_1979_ParticularSolutionsOfSGEIn2add1D}, the B\"acklund  transformation \cite{ChristiansenPL&OlsenOH_1979_OnDynamical2DSolutionsToSGE} and Painlev\'e transcendents \cite{KaliappanP&LakshmananM_1979_KadomtsevPetviashviliAnd2DSGEReductionToPainlevTranscendents}, etc. Moreover, numerical solutions for the $(2+1)$-dimensional undamped SGE have been proposed by \citet{ChristiansenPL&LomdahlPS_1981_NumercalSolutionOf2add1DimensionalSGSolitons} using a generalized leapfrog method, \citet{GuoBY&RodriguezPJ&VazquezL_1986_NumericalSolutionOfSGE} using two finite difference schemes, \citet{ArgyrisJ&HaaseM&HeinrichJC_1991_FiniteElementApproximationTo2DSGSolitons} using finite elements. \citet{XinJX_2000_ModelingLightBulletsWith2DSGE} studied the SGE as an asymptotic reduction of the two level dissipationless Maxwell-Bloch system, \citet{ShengQ&KhaliqAQM&VossDA_2005_NumericalSimulationOf2DSGSolitonsViaASplitCosineScheme} presented a numerical method with a split cosine scheme, and \citet{BratsosAG_2005_AnExplicitNumericalSchemeForSGEIn2add1D} used a three-time level fourth-order explicit finite difference scheme to solve the undamped SGE. Following a similar approach, \citet{BratsosAG_2007_TheSolutionOf2DSGEUsingTheMethodOfLines} transformed the SGE to a second-order initial value problem with the help of the method of lines. Numerical approaches for the damped SGE were proposed by \citet{NakajimaK&OnoderaYetal_1974_NumericalAnalysisOfVortexMotionOnJosephsonStructures} who considered dimensionless loss factors and unitless normalized bias, and \citet{GorriaC&GaidideiYBetal_2004_KinkPropagationAndTrappingIn2DCurvedJosephsonJunction} investigated nonlinear wave propagation in a planar wave guide consisting of two rectangular regions joined by a bent domain of constant curvature using as a model of the kink solution to the SGE. Additionally, \citet{DehghanM&MirzaeiD_2008_TheDualReciprocityBoundaryElementMethodFor2DSGE} developed the dual reciprocity boundary element method for both the undamped and damped $(2+1)$-dimensional SGE, and \citet{JiwariR&PanditS&MittalRC_2012_NumericalSimulationOf2DSGEByDifferentialQuadratureMethod} obtained a numerical scheme based on a polynomial differential quadrature method.

Although the SGE is nonintegrable except for $n=1$, some properties and exact solutions for the $(n+1)$-dimensional SGE have been obtained by various methods. \citet{KobayashiKK&IzutsuM_1976_ExactSolutionOnNDSGE} extensively studied the exact traveling wave solutions of the SGE in the field of theoretical physics. Many additional mathematical methods have been proposed for finding traveling wave solutions of the SGE. \citet{FengZS_2004_AnApproximateSGEandItsTravelingWaveSolutionInNadd1DSpace} applied the Painlev\'{e} analysis to the study of an approximate SGE and its traveling solitary wave solution in $(n+1)$-dimensional space. With the help of exact solutions to the cubic nonlinear Klein-Gordon fields, \citet{LouSY&HuHC&TangXY_2005_InteractionsAmongPeriodicWavesAndSolitaryWavesOfNadd1DSGField} studied the exact solutions for the $(n+1)$-dimensional SGE. Adopting $\beta = 0$ and $\psi(\boldsymbol{x}) = 1$, \citet{HozF&VadilloF_2012_NumericalSimulationOfNDSGEViaOperationalMatrices} generalized the exact soliton solution for the $(n+1)$-dimensional SGE:
\begin{equation}\label{eq:hdmrhc.exactsge}
 u(\boldsymbol{x},t)= 4\arctan\left[C\exp\left(\sum_{i=1}^n a_i x_i - b t\right)\right]
\end{equation}
where $\sum_{i=1}^n a_i^2=1+ b^2$. By adopting the proper \emph{ansatz}, more general solutions can be obtained for the multi-dimensional SGE, including the three-dimensional case allowing for non-constant $C$ \cite{AeroEL&BulyginAN&PavlovYV_2009_SolutionsOf3DSGE}. Obtaining the exact solution for the general SGE would be ideal, but unfortunately it is very difficult for practical engineering problems that are usually complex in nature. Despite numerical methods commonly  used in many types of linear and nonlinear PDEs, \citet{HozF&VadilloF_2012_NumericalSimulationOfNDSGEViaOperationalMatrices} remarked that there were no references to the numerical treatment of the SGE for dimensions larger than three, which motivated them to propose a numerical method for the $n$-dimensional SGE based on using operational matrices.

Many standard numerical methods for solving PDEs are widely used in engineering, but they usually require the construction and update of a mesh, which is an inherent disadvantage. In order to overcome these difficulties, recently the meshless numerical method has attracted attention. This method can establish a system of algebraic equations over the entire problem domain without using a predefined mesh. Rather, a set of scattered nodes, called field nodes, are used within the problem domain as well as on the boundaries of the domain \cite{LiuGR&GuYT_2005_AnIntroductionToMeshfreeMethodsAndTheirPrograming}. The meshless method does not require \emph{a priori} information about the relationship between the nodes for the interpolation or approximation of unknown functions over the field of variables \cite{LiuGR&GuYT_2005_AnIntroductionToMeshfreeMethodsAndTheirPrograming}. Employing the meshless method, \citet{DehghanM&ShokriA_2008_ANumericalMethodFor1DNonlinearSGECollocationAndRBF} studied the one-dimensional nonlinear SGE and used Thin Plate Spline Radial Basis Functions (TPS-RBF) to approximate the solution, and they also applied the TPS-RBF method to both the Klein-Gordon equation \cite{DehghanM&ShokriA_2009_NumericalSolutionOfNonlinearKGEUsingRBF} and the two-dimensional SGE \cite{DehghanM&ShokriA_2008_ANumericalMethodFor2DSGEUsingRBF}. A series of meshless approaches have been presented \cite{AsgariZ&HosseiniSM_2013_NumericalSolutionOf2DSGandMBEModelsUsingFourierSpectral, DehghanM&GhesmatiA_2010_NumericalSimulationOf2DSGSolitonsViaLocalWeakMeshlessTechniqueBasedOnRPIM, JiangZW&WangRH_2012_NumericaSolutionOf1DSGEUsingHighAccuracyMQQuasiInterpolation, JiwariR&PanditS&MittalRC_2012_ADifferentialQuadratureAlgorithmToSolve2DLinearHyperbolicTelegraphEquationWithDirichletAndNeumannBoundaryConditions, KaramanliA&MuganA_2013_StrongFormMeshlessImplementationOfTaylorSeriesMethod, MirzaeiD&DehghanM_2009_ImplementationOfMeshlessLBIEMethodTo2DNonlinearSGproblem, MirzaeiD&DehghanM_2010_MeshlessLocalPetrovGalerkinApproximationTo2DSGE, PekmenB&TezerSezginM_2012_DifferentialQuadratureSolutionOfNonlinearKGandSGE, ShaoWT&WuXH_2014_NumericalSolutionOfNonlinearKGandSGEUsingChebyshevTauMeshlessMethod, TaleeiA&DehghanM_2014_APseudoSpectralMethodThatUsesAnOverlappingMultidomainTechniqueForNumericalSolutionOfSGEIn1and2D}.
Moreover, since the nodal distribution for most existing meshless methods is preassigned, \citet{XuX_2015_cpc2015meshlessmethod} proposed a numerical two-step meshless method for soliton-like structures based on the optimal sampling density of kernel interpolation.

In dealing with high-dimensional PDEs, obtaining good quality approximate solutions is a difficult problem because of the so-called `curse of dimensionality'. High dimensional model representation (HDMR) \cite{LiGY_2001_HDMR, RabitzH&AlisOF_1999_GeneralFoundationsOfHighDimensionalModelRepresentations, RabitzH&AlisOF_2000_ManagingTyrannyOfParametersInMathematicalModelingOfPhysicalSystems, RabitzH&AlisOF_1998_EfficientInputOutputModelRepresentations, SobolIM_1993_SensitivityEstimatesForNonlinearMathematicalModels} provides a viable approach based on the fact that high-dimensional functions often can be efficiently expressed as sums of low-dimensional functions. The HDMR decomposition is also well known in statistics as the ANOVA (analysis of variance) decomposition \cite{EfronB_1981_anova, FisherRA_1925_anova, GriebelM_2010_anova, StoneC_1994_anova}. In recent years, the HDMR decomposition has been under rapid development becoming an important tool for understanding high-dimensional functions \cite{GriebelM_2010_anova, GriebelMetal_2010_anova, GriebelMetal_2013_smoothingANOVA, LiGY_2001_HDMR, LuoXP&LuZZ&XuX_2014_NonparametricKernelEstimationForANOVADecompositionAndSensitivityAnalysis, LuoXP&LuZZ&XuX_2014_ReproducingKernelTechniqueForHighDimensionalModelRepresentations, RabitzH&AlisOF_1999_GeneralFoundationsOfHighDimensionalModelRepresentations, RabitzH&AlisOF_2000_ManagingTyrannyOfParametersInMathematicalModelingOfPhysicalSystems, RabitzH&AlisOF_1998_EfficientInputOutputModelRepresentations, SobolIM_1993_SensitivityEstimatesForNonlinearMathematicalModels}. In this paper, we will use a HDMR decomposition in conjunction with the Fourier hyperbolic cross (HC) approximation \cite{LuoXP_2016_fourierhdmrhc}.

The remainder of the paper is organized as follows: Section \ref{sec:hdmrhc.funapp} presents a function approximation method using HDMR-HC. Then a new meshless numerical scheme is proposed in Section \ref{sec:hdmrhc.solsge} for solving the $(n+1)$-dimensional SGE using the HDMR-HC approximation. In Section \ref{sec:hdmrhc.eg}, we provide several examples with a comparative numerical error analysis. Section \ref{sec:hdmrhc.conclusion} summarizes the relevant results.

\section{Function approximation using HDMR-HC and the partition of unity}
\label{sec:hdmrhc.funapp}

\subsection{HDMR-HC approximation}
\label{sec:hdmrhc.funapp.app}

Let $n\geqslant 2$, $\lambda\in\mathbb{R}_+$ and $f(\boldsymbol{x})\in\mathcal{W}_1^s(\mathbb{T}^n_\lambda), s\in\mathbb{N}_0$ be a $n$-variate function which is $\frac{1}{\lambda}$-periodic in each variable, where $\mathbb{T}^n_\lambda$ is the $n$-torus given by \cite{LuoXP_2016_fourierhdmrhc}
\begin{equation*}
    \mathbb{T}^n_\lambda:=\left[-\frac{1}{2\lambda},\frac{1}{2\lambda}\right)^n
\end{equation*}
and the function space $\mathcal{W}_1^s(\mathbb{T}^n_\lambda)$ is defined by
\begin{equation}\label{eq:W1s}
    \mathcal{W}_1^s(\mathbb{T}^n_\lambda)
    =\{f\in\mathcal{L}(\mathbb{T}^n_\lambda):
    \partial^rf\in\mathcal{L}(\mathbb{T}^n_\lambda), ~|r|_\infty\leqslant s\}
\end{equation}
with the norm
\begin{equation*}
    \|f\|_{\mathcal{W}_1^s(\mathbb{T}^n_\lambda)}
    =\sum_{0\leqslant|r|_\infty\leqslant s}
    \|\partial^rf\|_{\mathcal{L}(\mathbb{T}^n_\lambda)}
\end{equation*}
where $r\in\mathbb{N}_0^n$ denotes a $n$-dimensional multi-index with the norm
\begin{equation}\label{eq:miN}
    |r|_\infty:=\max_{1\leqslant j\leqslant n} r_j
\end{equation}
and the $r$-th order mixed derivative is given by
\begin{equation}\label{eq:mixD}
    \partial^rf:=\frac{\partial^{|r|_1}f}{\partial x_1^{r_1}\cdots\partial x_n^{r_n}}
\end{equation}
We consider the multivariate Fourier series of $f(\boldsymbol{x})$
\begin{equation}\label{eq:mfs}
    f(\boldsymbol{x})\sim\sum_{m\in\mathbb{Z}^n}c_m \mathrm{e}^{2\pi\mathrm{i}\lambda m\cdot \boldsymbol{x}}, ~\boldsymbol{x}\in\mathbb{T}_\lambda^n
\end{equation}
where the Fourier coefficients $c_m$ are defined by
\begin{equation*}
    c_m=c_m(f)=\lambda^n\int_{\mathbb{T}_\lambda^n}
    f(\boldsymbol{x})\mathrm{e}^{-2\pi \mathrm{i}\lambda m\cdot \boldsymbol{x}}\ud \boldsymbol{x}, ~m\in\mathbb{Z}^n.
\end{equation*}

For a nonempty set $\{j_1,\cdots,j_v\}\subset\{1,\cdots,n\}$, let
\begin{equation*}
 \begin{aligned}
    \Lambda_{j_1,\cdots,j_v}:= & \left\{ m\in\mathbb{Z}^n: m_l\neq0~\forall~l\in\{j_1,\cdots,j_v\}\right.\\
    &\left.~\textrm{and}~m_l=0~\forall~l\notin\{j_1,\cdots,j_v\},~l\in\{1,\cdots,n\}\right\}
 \end{aligned}
\end{equation*}
and
\begin{equation*}
    \Lambda_{j_1,\cdots,j_v}^{k_\lambda}:=
    \{m\in\Lambda_{j_1,\cdots,j_v}:k_\lambda\leqslant(2\pi)^v
    \lambda_{j_1}\cdots\lambda_{j_v}|m_{j_1}\cdots m_{j_v}|<k_\lambda+1,~k_\lambda\in\mathbb{N}\}.
\end{equation*}
Note that $\Lambda_0$ is the set consisting of the $n$-dimensional zero vector. Then $\mathbb{Z}^n$ can be decomposed into the following form

\begin{align}\label{eq:decomZd}
    \mathbb{Z}^n=\Lambda_0+\sum_{j}\Lambda_j+\sum_{j_1<j_2}
    \Lambda_{j_1,j_2}+\cdots+\Lambda_{j_1,\cdots,j_n}
\end{align}

and we refer to this as a HDMR decomposition of $\mathbb{Z}^n$; then a multiple Fourier series can be decomposed with an HDMR structure
\begin{equation}\label{eq:mfsD}
 \begin{aligned}
     \sum_{m\in\mathbb{Z}^n}c_m \mathrm{e}^{2\pi i\lambda m\cdot x}
    = &~c_0+\sum_{k_\lambda=1}^\infty\left(\sum_{j}\!\sum_{m\in\Lambda_j^{k_\lambda}}c_m\mathrm{e}^{2\pi i\lambda m\cdot x}\right.\\
    & \left.+\!\!\sum_{j_1<j_2}\!\sum_{m\in\Lambda_{j_1,j_2}^{k_\lambda}}\!\!\!\!\!\! c_m\mathrm{e}^{2\pi i\lambda m\cdot x}+\cdots+\hspace{-5mm} \sum_{m\in\Lambda_{j_1,\cdots,j_n}^{k_\lambda}}\hspace{-5mm}
    c_m\mathrm{e}^{2\pi i\lambda m\cdot x}\right)
 \end{aligned}
\end{equation}
If let
\begin{equation}\label{eq:Lambda}
    \Lambda^{K_\lambda}=\Lambda_0
    +\sum_{k_\lambda=1}^{K_\lambda}\left(\sum_{j}\Lambda_j^{k_\lambda}
    +\sum_{j_1<j_2}\Lambda_{j_1,j_2}^{k_\lambda}+\cdots+
    \Lambda_{j_1,\cdots,j_d}^{k_\lambda}\right),
\end{equation}
that is,
\begin{equation}\label{eq:hdmrhc.LK}
    \Lambda^{K_\lambda}:=\left\{m\in\mathbb{Z}^n:\prod_{j=1}^n
    \max\{2\pi\lambda|m_j|,1\}\leqslant K_\lambda\right\},
\end{equation}
then we define the Fourier HDMR-HC partial sum up to $K_\lambda$-th order \cite{LuoXP_2016_fourierhdmrhc}
\begin{equation}\label{eq:HDMR-HC}
    S_{K_\lambda}(\boldsymbol{x},f)=\sum_{m\in\Lambda^{K_\lambda}}
    c_m\mathrm{e}^{2\pi \mathrm{i}\lambda m\cdot \boldsymbol{x}}, ~\boldsymbol{x}\in\mathbb{T}_\lambda^n.
\end{equation}
Suppose $M$ is the number of points $\boldsymbol{x}\in\mathbb{T}_\lambda^n$, which depends on both $\lambda$ and $n$, and
\begin{displaymath}
 F_{\lambda,m}(\boldsymbol{x}) = \mathrm{e}^{2\pi\mathrm{i}\lambda m\cdot\boldsymbol{x}}
\end{displaymath}
then Eq.\eqref{eq:HDMR-HC} can be rewritten as
\begin{equation}\label{eq:HDMR-HC2}
    S_{K_\lambda}(\boldsymbol{x},f)=\sum_{m=1}^{M}c_m F_{\lambda,m}(\boldsymbol{x}), ~\boldsymbol{x}\in\mathbb{T}_\lambda^n
\end{equation}
Moreover, it follows that, for any $f\in\mathcal{W}_1^s(\mathbb{T}^n_\lambda)$, $s>\gamma+p+1$ and $K_\lambda\geqslant\delta$, the bound
\begin{equation*}
    |\partial^rS_{K_\lambda}(\boldsymbol{x},f)-\partial^rf(\boldsymbol{x})|\leqslant C_{n,\lambda,s,p,\gamma,\delta}
    \|f\|_{\mathcal{W}_1^s(\mathbb{T}^n_\lambda)}K_\lambda^{-(s-\gamma-p-1)}, ~~|r|_\infty=p
\end{equation*}
holds pointwise almost everywhere, where the constant $C_{n,\lambda,s,p,\gamma,\delta}$ depends on $n,\lambda,s,p$ and $\delta$.

From the definition of $\mathbb{T}^n_\lambda$, we find that $\mathbb{T}^n_\lambda$ is smaller as $\lambda$ becomes larger. The convergence rate of the Fourier HDMR-HC partial sums $S_{K_\lambda}$ is closely related to the value of $\lambda$. In particular, for a given accuracy, when $\lambda$ is large, $S_{K_\lambda}$ generally is well approximated by a low order truncated HDMR, and this is the basic starting point of the Fourier HDMR-HC approximation. For example, if there is a function of $n=5$ variables, and suppose $\lambda=1/\pi$ and the desired accuracy is $\varepsilon=\mathbf{O}(K_\lambda^{-(s-p-\gamma-1)})$, then a $K_\lambda$-th order Fourier HDMR-HC partial sum $S_{K_\lambda}$ is just a $2$nd order truncated HDMR of the function when $K_\lambda\leqslant31$. Therefore, we expect that a low order truncated HDMR-HC can be used to effectively capture the behavior of a high-dimensional function and its derivatives.

\subsection{Partition of unity}
\label{subsec:hdmrhc.funapp.pu}
In this subsection, we will discuss how to approximate a function via HDMR-HC. Suppose $u(\boldsymbol{x})$ is a function defined  on $\mathbb{R}$, where $\boldsymbol{x} \in \Omega\subseteq\mathbb{R}^n$. First, we divide the domain $\Omega$ into $D \in \mathbb{N}$ subdomains, denoted as $\Omega_j,~(j=1,2,\cdots,D)$, thus $\Omega =\bigcup_{j=1}^D\Omega_j$. In each subdomain $\Omega_j$, let $\boldsymbol{x}_j$ be the centre of the region, and let $\boldsymbol{\chi}_j$ denote all $M_j$ nodes in $\Omega_j$, i.e., $\sum_{j=1}^DM_j=N$. Further, suppose $u(\boldsymbol{x}),\boldsymbol{x} \in \Omega_j$ can be represented as $u_j$, then we have
\begin{equation}
\label{eq:hdmrhc.fj}
 u_j(\boldsymbol{x}) = J_\lambda(\boldsymbol{x}-\boldsymbol{x}_j)u(\boldsymbol{x}) = J_j u(\boldsymbol{x}),~ \boldsymbol{x}\in \Omega_j
\end{equation}
where
\begin{equation}
  \Omega_j\subset\textrm{supp}~J_\lambda(\boldsymbol{x}-\boldsymbol{x}_j)~~~\textrm{and}~~~
  \textrm{supp}~J_\lambda(\boldsymbol{x})\subset\mathbb{T}^n_\lambda
\end{equation}
is the characteristic function of $u_j$ satisfying
\begin{equation}
\label{eq:hdmrhc.cfuncnd}
 \sum_{j=1}^{D} J_j=1,~\forall \boldsymbol{x} \in\Omega\subseteq\mathbb{R}^n
\end{equation}
From Eq.\eqref{eq:HDMR-HC2}, $u_j$ at any point $\boldsymbol{x} \in \Omega_j$ can be approximated as
\begin{equation}
\label{eq:hdmrhc.appfj}
 \hat{u}_j(\boldsymbol{x}) = \sum_{m_j=1}^{M_j}c_{m_j} F_{\lambda,m_j}(\boldsymbol{x}-\boldsymbol{x}_j) =\boldsymbol{F_{\lambda}}_j^\mathrm{T}(\boldsymbol{x})\boldsymbol{c}_{j}
\end{equation}
where $\boldsymbol{F_{\lambda}}_j^\mathrm{T}(\boldsymbol{x})=[F_{\lambda,1}, F_{\lambda,2}, \cdots, F_{\lambda,M_j}]$ with the same $\lambda$ for all nodes and $c_{m_j}$ is the $m_j$-th unknown coefficient. Utilizing the values of these $M_j$ nodes, there are $M_j$ equations with one for each node, then we have following matrix form
\begin{equation}
\label{eq:hdmrhc.appfimatrix}
 \boldsymbol{U}_j = \boldsymbol{F}_j\boldsymbol{c}_{j}
\end{equation}
where
\begin{displaymath}
 \boldsymbol{U}_j = [u_{1}, u_{2},\cdots,u_{M_j}]^\mathrm{T}
\end{displaymath}
is the vector of function values at the $M_j$ nodes, and
\begin{displaymath}
 \boldsymbol{c}_{j} = [c_{1}, c_{2},\cdots,c_{M_j}]^\mathrm{T}
\end{displaymath}
is the vector of undetermined coefficients with
\begin{displaymath}
 \boldsymbol{F}_j =
 \left[ \begin{array}{cccc}
  F_{\lambda,11} & F_{\lambda,12} & \cdots & F_{\lambda,1M_j} \\
  F_{\lambda,21} & F_{\lambda,22} & \cdots & F_{\lambda,2M_j} \\
  \vdots & \vdots &  & \vdots  \\
  F_{\lambda,M_j1} & F_{\lambda,M_j2} & \cdots & F_{\lambda,M_jM_j}
 \end{array} \right]
\end{displaymath}
where $F_{\lambda,im_j} = F_{\lambda,m_j}(\boldsymbol{\chi}_i)$. Suppose $\boldsymbol{F}_j^{-1}$ exists (i.e., this condition can always be satisfied \cite{PowellMJD_1992_TheoryOfRBFApproximationIn1990, WendlandH_1998_ErrorEstimatesForInterpolationByCompactlySupportedRBFOfMinimalDegree}), then $\boldsymbol{c}_{j}$ can be obtained by solving Eq.\eqref{eq:hdmrhc.appfimatrix}, i.e.
\begin{equation}\label{eq:hdmrhc.hdmrcoej}
  \boldsymbol{c}_{j} = \boldsymbol{F}_j^{-1}\boldsymbol{U}_j
\end{equation}
From Eq.\eqref{eq:hdmrhc.fj}, we have
\begin{equation}\label{eq:hdmrhc.uuj}
 \boldsymbol{U}_j = J_j \boldsymbol{U}
\end{equation}
where $\boldsymbol{U}= \bigcup_{j=1}^{D}\boldsymbol{U}_j$. Then
\begin{equation}\label{eq:hdmrhc.hdmrcoej2}
  \boldsymbol{c}_{j} = \boldsymbol{F}_j^{-1}J_j \boldsymbol{U}
\end{equation}

Substitute both Eq.\eqref{eq:hdmrhc.hdmrcoej} and Eq.\eqref{eq:hdmrhc.uuj} back into Eq.\eqref{eq:hdmrhc.appfj}, we have an approximation of function $u_j$:
 \begin{equation}\label{eq:hdmrhc.kappp}
 \hat{u}_j(\boldsymbol{x}) = \boldsymbol{F_\lambda}_j^\mathrm{T}(\boldsymbol{x})\boldsymbol{F}_j^{-1}J_j \boldsymbol{U}
\end{equation}

Furthermore, from Eqs.\eqref{eq:hdmrhc.fj} and \eqref{eq:hdmrhc.cfuncnd} we have
 \begin{equation}
 u(\boldsymbol{x}) = \sum_{j=1}^{D} J_j u(\boldsymbol{x}) = \sum_{j=1}^{D} u_j(\boldsymbol{x}), ~\boldsymbol{x}\in\Omega
\end{equation}
Then $u(\boldsymbol{x})$ at any point $\boldsymbol{x}\in\Omega$ can be formally approximated as
\begin{equation}\label{eq:hdmehc.fmhdmrapp}
 \begin{aligned}
  \hat{u}(\boldsymbol{x}) &= \sum_{j=1}^{D} \boldsymbol{F_\lambda}_j^\mathrm{T}(\boldsymbol{x}) \boldsymbol{F}_j^{-1}J_j \boldsymbol{U}\\
  &= \sum_{j=1}^{D} J_j \boldsymbol{F_\lambda}_j^\mathrm{T}(\boldsymbol{x}) \boldsymbol{F}_j^{-1} \boldsymbol{U}\\
  &= \boldsymbol{\Phi}^\mathrm{T}(\boldsymbol{x})\boldsymbol{U}
 \end{aligned}
\end{equation}

\section{Solution for $(n+1)$-dimensional SGEs}
\label{sec:hdmrhc.solsge}

Now, we present the numerical scheme for solving the $(n+1)$-dimensional SGE based on using the Fourier HDMR-HC approximation. Suppose the approximated function of the field function $u(\boldsymbol{x})$, ($\boldsymbol{x} = \{x_1,x_2,\cdots,x_n\}$) is formally denoted as
\begin{equation}\label{eq:hdmrhc.hdmrhcapp}
 \hat{u}(\boldsymbol{x}) = \sum_{m=1}^{M_j}c_{m_j} F_{\lambda,m_j}(\boldsymbol{x}-\boldsymbol{x}_j) =\boldsymbol{F_\lambda}_j^\mathrm{T}(\boldsymbol{x})\boldsymbol{F}_j^{-1}\boldsymbol{U_j}
\end{equation}
where $M_j$ is the number of field nodes used in the selected domain, and $\boldsymbol{U_j}$ is the vector that collects the true nodal function values for these $M_j$ field nodes, and $\boldsymbol{x}_j$ is the centre of this selected region. Further, the derivatives of $u(\boldsymbol{x})$ at any point $\boldsymbol{x}$ can be approximated as
\begin{equation}\label{eq:hdmrhc.derihdmrhcapp}
 \frac{\partial^{p} }{\partial x_l^p}\hat{u}(\boldsymbol{x}) = \frac{\partial^{p} \boldsymbol{F_\lambda}_j^\mathrm{T}(\boldsymbol{x})}{\partial x_l^p}\boldsymbol{F}_j^{-1}\boldsymbol{U_j}
\end{equation}
where $x_l$ denotes one element of $\boldsymbol{x} = \{x_1,x_2,\cdots,x_n\}$.

In this paper, the time derivatives are approximated by the time-stepping method and we have the following approximation:
\begin{equation}\label{eq:hdmrhc.timed2dis}
 \frac{\partial ^2 u}{\partial t^2} \approx \frac{1}{\tau^2} \left[ u^{(k+1)} - 2u^{(k)} + u^{(k-1)} \right]
\end{equation}
\begin{equation}\label{eq:hdmrhc.timed1dis}
 \frac{\partial u}{\partial t} \approx \frac{1}{2\tau} \left[ u^{(k+1)} - u^{(k-1)} \right]
\end{equation}
where $\tau$ is the time step, and $u^{(k)}$ is the approximate value of $u(\boldsymbol{x},t)$ at $(\boldsymbol{x},t_k)$, $t_k = k\tau$. Moreover the Crank-Nicolson scheme is used to approximate $u$ at three respective times as
\begin{equation}\label{eq:hdmrhc.timedCN}
 u(\boldsymbol{x},t) \approx \frac{1}{3}\left[u^{(k+1)} +u^{(k)} + u^{(k-1)}\right]
\end{equation}
To manage the nonlinearity, a Quasilinearization Method (QLM) is adopted. The QLM is very effective for dealing with the nonlinear aspects of the SGE and other PDEs. In this fashion the nonlinear term in Eq. \eqref{eq:hdmrhc.sge} can be represented as
\begin{equation}\label{eq:hdmrhc.qlmsge}
 \sin(u)= \sin(u^{(k)})+(u^{(k+1)}-u^{(k)})\cos(u^{(k)})
\end{equation}
Thus, Eq.\eqref{eq:hdmrhc.sge} can be discretized as
\begin{equation}\label{eq:hdmrhc.sgedis}
\begin{aligned}
 & \frac{1}{3}\Delta u^{(k+1)}-(\eta+\mu+\varphi^{(k)}) u^{(k+1)} \\
 =&-\frac{1}{3}\Delta u^{(k)}-(2\mu+\varphi^{(k)})u^{(k)}+\psi(\boldsymbol{x})\sin(u^{(k)})\\
  &-\frac{1}{3}\Delta u^{(k-1)}-(\eta-\mu) u^{(k-1)}
\end{aligned}
\end{equation}
where $\mu =\tau^{-2}$, $\eta =\frac{\beta}{2\tau}$ and $\varphi^{(k)}=\psi(\boldsymbol{x})\cos(u^{(k)})$.

Suppose $N$ field nodes are denoted as $\boldsymbol{\chi}=\{\boldsymbol{\chi}\}_1^{N}$, where $N$ is determined by both $D$ and $M_j$. Then from Eq.\eqref{eq:hdmehc.fmhdmrapp} we have an approximation of the field variable $u$ according to the HDMR-HC approximation
\begin{equation}\label{eq:hdmrhc.uappp}
  \hat{u}(\boldsymbol{x}) = \boldsymbol{\Phi}^\mathrm{T}(\boldsymbol{x})\boldsymbol{U}=\sum_{j=1}^{N}\phi_{j}(\boldsymbol{x})u_{j}
\end{equation}
where $\phi_{j}(\boldsymbol{x})$ depends on $J_j, \boldsymbol{F_\lambda}_j^\mathrm{T}(\boldsymbol{x})$ and $\boldsymbol{F}_j^{-1}$. The derivatives of $\hat{u}$ can be approximated as
\begin{equation}\label{eq:hdmrhc.duapp}
 \Delta \hat{u}(\boldsymbol{x}) =\Delta \boldsymbol{\Phi}^\mathrm{T}(\boldsymbol{x})\boldsymbol{U}
\end{equation}
Therefore, for any point $\boldsymbol{x}_i $, the approximation in Eq.\eqref{eq:hdmrhc.sgedis} can be written as
\begin{equation}\label{eq:hdmrhc.sgedisb}
\begin{aligned}
 & \sum_{j=1}^{N}\left[\frac{1}{3}\Delta \phi_{j}(\boldsymbol{x}_i)-(\eta+\mu+\varphi^{(k)})\phi_{j}(\boldsymbol{x}_i)\right]u_{j}^{(k+1)} \\
 =&\sum_{j=1}^{N}\left[-\frac{1}{3}\Delta \phi_{j}(\boldsymbol{x}_i)-(2\mu+\varphi^{(k)}) \phi_{j}(\boldsymbol{x}_i)\right]u_{j}^{(k)} \\
 &~+\sum_{j=1}^{N}\left[-\frac{1}{3}\Delta \phi_{j}(\boldsymbol{x}_i)-(\eta-\mu)\phi_{j}(\boldsymbol{x}_i)\right]u_{j}^{(k-1)}\\
 &~+\psi(\boldsymbol{x}_i)\sin(u^{(k)})
\end{aligned}
\end{equation}
Let $A_{ij} = \frac{1}{3}\Delta \phi_{j}(\boldsymbol{x}_i)$, $B_{ij} = \phi_{j}(\boldsymbol{x}_i) $, and $E_{ij}^{(k)} = A_{ij} - (\eta+\mu+\varphi^{(k)})B_{ij}$, $G_{ij}^{(k)} = -A_{ij} -(2\mu+\varphi^{(k)}) B_{ij}$, $H_{ij} = A_{ij} - (\eta-\mu)B_{ij}$, $C_i^{(k)} = \psi(\boldsymbol{x}_i)\sin(u^{(k)})$, then Eq.\eqref{eq:hdmrhc.sgedisb} can be re-written as
\begin{equation}\label{eq:hdmrhc.sgedisc}
 \sum_{j=1}^{N}E_{ij}^{(k)}u_{j}^{(k+1)}=\sum_{j=1}^{N}G_{ij}^{(k)}u_{j}^{(k)} + \sum_{j=1}^{N}H_{ij}u_{j}^{(k-1)}+ C_i^{(k)}
\end{equation}
For all $N$ field nodes $\boldsymbol{\chi}$ we have following matrix form:
\begin{equation}\label{eq:hdmrhc.sgematrix}
 \boldsymbol{E}^{'(k)}\boldsymbol{\hat{u}}^{'(k+1)}=\boldsymbol{G}^{'(k)}\boldsymbol{\hat{u}}^{'(k)}+ \boldsymbol{H}'\boldsymbol{\hat{u}}^{'(k-1)}+ \boldsymbol{C}^{'(k)}
\end{equation}
where $\boldsymbol{C}^{'(k)}= [C_1^{(k)},C_2^{(k)},\cdots,C_N^{(k)}]^\mathrm{T}$.

In using a meshless strong method to solve the PDE, the solution can be unstable if there is a derivative boundary condition, so the fictitious points method is used to impose derivative boundary conditions \cite{LiuGR&GuYT_2005_AnIntroductionToMeshfreeMethodsAndTheirPrograming}. Suppose there are $N_b$ nodes on the boundary, then along the derivative boundaries, another $N_b$ fictitious points are added outside of the domain. Two sets of equations are established at each derivative boundary node: one for the derivative boundary condition, and the other for the governing equation. With the $N_b$ additional degrees of freedom, $\{u_{(N+1)},u_{(N+2)},\cdots,u_{(N+N_b)}\}$, added into the system, then Eq.\eqref{eq:hdmrhc.sgedisb} can be re-written as
\begin{equation}\label{eq:hdmrhc.sgediscbc1}
 \sum_{j=1}^{N+N_b}E_{ij}^{(k)}u_{j}^{(k+1)}=\sum_{j=1}^{N+N_b}G_{ij}^{(k)}u_{j}^{(k)} +\sum_{j=1}^{N+N_b}H_{ij}u_{j}^{(k-1)}+ C_i^{(k)}
\end{equation}
and for a node at $\boldsymbol{x}_{i_b}$, that is, on the derivative boundary, the derivative boundary conditions have the form
\begin{equation}\label{eq:hdmrhc.sgediscbc2}
 \sum_{j=1}^{N+N_b}\!\!\frac{\partial}{\partial l} \phi_{j}(\boldsymbol{x}_{i_b}) u_{j}^{(k+1)}= -\!\!\sum_{j=1}^{N+N_b}\!\!\frac{\partial}{\partial l} \phi_{j}(\boldsymbol{x}_{i_b})u_{j}^{(k)}-\!\!\sum_{j=1}^{N+N_b}\!\!\frac{\partial}{\partial l} \phi_{j}(\boldsymbol{x}_{i_b})^{(k-1)}+ 3w(\boldsymbol{x}_{i_b},t)
\end{equation}
Assembling Eqs.\eqref{eq:hdmrhc.sgediscbc1} and \eqref{eq:hdmrhc.sgediscbc2} for the corresponding nodes, the discretized global system equation becomes
\begin{equation}\label{eq:hdmrhc.sgematrixbc}
 \boldsymbol{E}^{(k)}\boldsymbol{\hat{u}}^{(k+1)}=\boldsymbol{G}^{(k)}\boldsymbol{\hat{u}}^{(k)} +\boldsymbol{H}\boldsymbol{\hat{u}}^{(k-1)}+ \boldsymbol{C}^{(k)}
\end{equation}
where $\boldsymbol{E}$, $\boldsymbol{G}$, $\boldsymbol{B}$ and $\boldsymbol{H}$ are ${(N+N_b)\times(N+N_b)}$ matrices, $\boldsymbol{\hat{u}} = [\hat{u}_{1}, \hat{u}_{2}, \cdots ,\hat{u}_{(N+N_b)}]^\mathrm{T}$, $\boldsymbol{C}^{(k)}=[C_1^{(k)},C_2^{(k)},\cdots,C_N^{(k)},3w(\boldsymbol{x}_{N+1},t),3w(\boldsymbol{x}_{N+2},t),\cdots,3w(\boldsymbol{x}_{N+N_b},t)]^\mathrm{T}$.
At the first time level, i.e. $k=0$, we adopt the following:
\begin{equation}\label{eq:hdmrhc.initial1}
 \boldsymbol{\hat{u}}^{(0)} = \boldsymbol{v_1}
\end{equation}
and
\begin{equation}\label{eq:hdmrhc.initial2}
 \boldsymbol{\hat{u}}^{(-1)} = \boldsymbol{\hat{u}}^{(1)}-2 \tau \boldsymbol{v_2}
\end{equation}
where $\boldsymbol{v_1}$ and $\boldsymbol{v_2}$ are the initial conditions for all nodes $\boldsymbol{\chi}$ introduced in Eqs.\eqref{eq:hdmrhc.sgeic1} and \eqref{eq:hdmrhc.sgeic2}.

\section{Numerical experiments}
\label{sec:hdmrhc.eg}

In this section, the proposed meshless numerical scheme is applied to several examples to show the efficiency and accuracy for the $(n+1)$-dimensional SGE. As mentioned in the previous section, to approximate the time derivatives we use a finite difference method, so an iterative scheme is employed to reach the final time $t$. In order to test the performance of the numerical solution, we use the $L_{\infty}$ error and root-mean-square (RMS) error norms defined as
\begin{equation}\label{eq:hdmrhc.linferr}
  L_{\infty} =\left\|f(x_i)-\hat{f}(x_i)\right\|_{\infty} = \max_{1\leqslant i \leqslant N}\left|f(x_i)-\hat{f}(x_i)\right|
\end{equation}
and
\begin{equation}\label{eq:hdmrhc.rmserr}
  \textrm{RMS} = \sqrt{\frac{1}{N}\sum_{i=1}^N\left|f(x_i)-\hat{f}(x_i)\right|^2}
\end{equation}
where $N$ is the number of nodes, $f(x_i)$ is the exact solution, and $\hat{f}(x_i)$ is the numerical solution. To assess both the stability and the solution accuracy, we compute the condition number of the system matrix, which is defined as
\begin{equation}\label{eq:sgker.cndn}
 \kappa(\boldsymbol{F})=\|\boldsymbol{F}\|\|\boldsymbol{F}^{-1}\|
\end{equation}
where $\kappa(\boldsymbol{F})$ depends on the parameter $\lambda$ and the number of nodes $N$.

\subsection{Test problem for a $(2+1)$-dimensional SGE}
\label{subsec:hdmrhc.eg.test2D}

The test problem for a $(2+1)$-dimensional SGE has the following form \cite{DjidjeliK&PriceWG&TwizellEH_1995_NumericalSolutionsOfADampedSGEinTwoSpaceVariables, DehghanM&ShokriA_2008_ANumericalMethodFor2DSGEUsingRBF, JiwariR&PanditS&MittalRC_2012_NumericalSimulationOf2DSGEByDifferentialQuadratureMethod}
\begin{equation}\label{eq:hdmrhc.test2D}
 \begin{split}
  & u_{tt} = \Delta u - \sin(u), \\
  &~~(x,y) \in \Omega,~t>0,
  \end{split}
\end{equation}
where $\Omega=[-7,7]\times[-7,7]$, and the initial conditions are
\begin{equation}\label{eq:hdmrhc.test2Dic}
\begin{split}
  & v_1(x,y) = 4 \tan^{-1}\left[\exp(x+y)\right], ~(x,y) \in \Omega\\
  & v_2(x,y) = -\frac{4 \exp (x+y)}{1+\exp (2x+2y)},~(x,y) \in \Omega.
  \end{split}
\end{equation}
and the Neumann boundary condition is
\begin{equation}\label{eq:hdmrhc.test2Dbc}
 w(x,y,t) = \frac{4 \exp (x+y-t)}{1+\exp (2x+2y-2t)}, ~(x,y) \in \partial \Omega,~t>0
\end{equation}

The analytic solution of this problem is:
\begin{equation}\label{eq:hdmrhc.test2Dsol}
 u(x,y,t) = 4 \tan^{-1}\left[ \exp (x+y-t) \right],~(x,y) \in \Omega
\end{equation}

In this example, both the proposed HDMR-HC meshless method and the radial basis function (RBF) method in \cite{DehghanM&ShokriA_2008_ANumericalMethodFor2DSGEUsingRBF} are used to numerically solve the equation. Since the field nodes of the RBF method in \cite{DehghanM&ShokriA_2008_ANumericalMethodFor2DSGEUsingRBF} is the Sobol sequence with $N=3249$, then we adopt the same total number of field nodes for the HDMR-HC meshless method ($D=49$). The time step $\tau$ is set to $0.001$. The results of the two different measures of error are presented in Table \ref{tab:hdmrhc.test2Derr}. We see that with the same number of field nodes, the errors of the proposed HDMR-HC meshless method are smaller than those of the RBF method. The condition numbers at particular times are also listed in Table \ref{tab:hdmrhc.test2Derr}.

\begin{table}[!htp]
\small
\begin{threeparttable}
 \caption{Errors and condition number $\kappa$ for $(2+1)$-D SGE}
 \label{tab:hdmrhc.test2Derr}
 \begin{tabular}{cccccc}
  \toprule
    & \multicolumn{2}{c}{$L_\infty$-error} & \multicolumn{2}{c}{RMS-error} & \\
    \raisebox{1.5ex}[0cm][0cm]{Time($s$)}  & RBF & HDMR-HC & RBF & HDMR-HC &\raisebox{1.5ex}[0cm][0cm]{$\kappa$}\\
    \midrule
    1.0 & 0.0670 & 0.0326 & 0.0050 & 0.0043 & $3.4\times10^5$\\
    3.0 & 0.0834 & 0.0343 & 0.0103 & 0.0045 & $5.0\times10^5$ \\
    5.0 & 0.1015 & 0.0355 & 0.0145 & 0.0045 & $4.2\times10^5$ \\
    7.0 & 0.1516 & 0.0368 & 0.0187 & 0.0047 & $6.3\times10^5$ \\
  \bottomrule
 \end{tabular}
 \small Note: The results of RBF method come from Ref. \cite{DehghanM&ShokriA_2008_ANumericalMethodFor2DSGEUsingRBF}
\end{threeparttable}
\end{table}

\subsection{Test problem for a $(5+1)$-dimensional SGE}
\label{subsec:hdmrhc.eg.test5D}

To further test the proposed HDMR-HC scheme, we choose a $(5+1)$-dimensional example, which involves all the implementation issues explained in the previous subsection. The exact solution has the form of Eq.\eqref{eq:hdmrhc.exactsge}
\begin{displaymath}
 u(\boldsymbol{x},t)= 4\arctan\left[C\exp\left(\sum_{i=1}^5 a_i x_i - b t\right)\right], ~\boldsymbol{x}=(x_1,x_2,\cdots,x_5)\in \Omega, ~t>0
\end{displaymath}
where $C=1$, $b=1$, $\Omega=[-6,6]\times[-6,6]\times\cdots\times[-6,6]\subseteq \mathbb{R}^5$ and
\begin{equation}\label{eq:hdmrhc.testsol5D}
 a_i = \left\{\begin{array}{cl}
         \sqrt{2}/2, & i=1,2,3 \\[3pt]
         1/2, & i=4,5
       \end{array}\right.
\end{equation}
The initial conditions are
\begin{equation}\label{eq:hdmrhc.test5Dic}
\begin{aligned}
  & v_1(\boldsymbol{x}) = 4\arctan\left[C\exp\left(\sum_{i=1}^5 a_i x_i\right)\right],~\boldsymbol{x}\in \Omega,\\
  & v_2(\boldsymbol{x}) = -\frac{4bC \exp \left( \sum_{i=1}^n a_i x_i\right)}{1+C^2\exp \left( \sum_{i=1}^n 2a_i x_i\right)},~\boldsymbol{x}\in \Omega,
  \end{aligned}
\end{equation}
and the Neumann boundary conditions are
\begin{equation}\label{eq:hdmrhc.test5Dbc}
   w_j(\boldsymbol{x},t) = \frac{4a_jC \exp \left( \sum_{i=1}^5 a_i x_i- bt\right)}{1+C^2\exp \left( \sum_{i=1}^5 2a_i x_i- 2bt\right)},~j=1,2,\cdots,5,~\boldsymbol{x} \in \partial \Omega,~t>0
\end{equation}

Both the proposed HDMR-HC method and the RBF method are used to solve the equation. The Sobol sequence is chosen as the field nodes with a total number of $N=2^{16}$ for the HDMR-HC method ($D=243$ and we use the HDMR approximation up to order $3$) and $N'=2^{18}$ for the RBF method. In this case the time step is chosen as $\tau=0.2$. Table \ref{tab:hdmrhc.test5Dec} presents $L_\infty$, RMS errors and the condition number $\kappa$ at some selected times $t$.

\begin{table}[!htp]
 \caption{Errors and condition number $\kappa$ for $(5+1)$-D SGE}
 \label{tab:hdmrhc.test5Dec}
 \begin{tabular}{cccccc}
  \toprule
    & \multicolumn{2}{c}{$L_\infty$-error} & \multicolumn{2}{c}{RMS-error} &  \\
    \raisebox{1.5ex}[0cm][0cm]{Time($s$)} & RBF & HDMR-HC & RBF & HDMR-HC & \raisebox{1.5ex}[0cm][0cm]{$\kappa$} \\
    \midrule
    1.0 & 0.2071 & 0.1083 & 0.0130 & 0.0076 & $4.4\times10^6$ \\
    3.0 & 0.1956 & 0.0910 & 0.0143 & 0.0079 & $3.7\times10^6$ \\
    5.0 & 0.2132 & 0.1205 & 0.0165 & 0.0080 & $4.5\times10^6$ \\
    7.0 & 0.2203 & 0.1124 & 0.0187 & 0.0081 & $4.9\times10^6$ \\
  \bottomrule
 \end{tabular}
\end{table}

\section{Conclusions}
\label{sec:hdmrhc.conclusion}

In this paper, we propose a new meshless solution method for high-dimensional sine-Gordon equations. First, we present a function approximation using the HDMR-HC decomposition. Then we divided the whole domain into several subdomains with the help of the partition of unity, and obtain a function approximation at any random point in each subdomain. Hence, we develop a numerical procedure for the high-dimensional SGEs by a meshless strong solution method. The time-stepping method is used to approximate the time derivatives of SGEs, and a quasilinearization scheme is performed to treat the nonlinearity of the equation. Finally, to demonstrate the accuracy of the proposed method with two numerical experiments. The examples suggest that the proposed procedure is attractive for solving high-dimensional SGEs.

\begin{acknowledgements}
The authors X.X. and X.L. acknowledge support from the National Science Foundation (Grant No. CHE-1763198), and H.R. acknowledges support from the Templeton Foundation (Grant No. 52265).
\end{acknowledgements}

\bibliographystyle{spbasic} %_sortnum
\bibliography{referenceXu}

%\begin{thebibliography}{}
% etc
%\end{thebibliography}

\end{document}